\documentclass[fleqn]{mat01}
\usepackage{times,mathtimy,amssymb,epsfig}
\begin{document}

\makeatletter
\def\artpath#1{\def\@artpath{#1}}
\makeatother
\artpath{c:/prema/maths}

\setcounter{page}{79}
\firstpage{79}

\font\xxxx=tir at 9.4pt

\newtheorem{theo}{Theorem}
\renewcommand\thetheo{\arabic{section}.\arabic{theo}}
\newtheorem{theor}[theo]{\bf Theorem}
\newtheorem{lem}[theo]{Lemma}
\newtheorem{propo}{\rm PROPOSITION}
\newtheorem{rema}{Remark}
\newtheorem{defn}[theo]{\rm DEFINITION}
\newtheorem{exam}{Example}
\newtheorem{coro}[theo]{\rm COROLLARY}
\newtheorem{case}{Case}
\newtheorem{subcase}{Subcase}
\def\conjecture{\trivlist\item[\hskip\labelsep{\it Conjecture.}]}

\renewcommand{\theequation}{\thesection\arabic{equation}}

\title{Inverse solutions for a second-grade fluid for porous medium
channel and Hall current effects}

\markboth{Muhammad R~Mohyuddin and Ehsan Ellahi Ashraf}{Inverse
solutions for a second-grade fluid}

\author{MUHAMMAD R~MOHYUDDIN$^{1,2}$ and EHSAN ELLAHI ASHRAF$^{3}$}

\address{$^{1}$Corresponding author: Department of Mathematics,
Quaid-i-Azam University, 45320, Islamabad 44000, Pakistan\\
\noindent $^{2}$Present address: Department of Condensed Matter Physics,
Main Building, Strada~Costeira 11, 34014, ICTP, Trieste, Italy\\
\noindent $^{3}$College of Aeronautical Engineering, National University
of Sciences and Technology, PAF Academy, Risalpur 24090, Pakistan\\
\noindent E-mail: $^{1}$m\_raheel@yahoo.com; $^{2}$mmohyudd@ictp.trieste.it}

\volume{114}

\mon{February}

\parts{1}

\Date{MS received 29 August 2003; revised 19 December 2003}

\begin{abstract}
Assuming certain forms of the stream function inverse solutions of an
incompressible viscoelastic fluid for a porous medium channel in the
presence of Hall currents are obtained. Expressions for streamlines,
velocity components and pressure fields are described in each case and
are compared with the known viscous and second-grade cases.
\end{abstract}


\maketitle

\section{Introduction}

In recent years the theoretical study of MHD channel flows has been a
subject of great interest due to its widespread applications in
designing cooling systems with liquid metals, petroleum industry,
purification of crude oil, polymer technology, centrifugal separation of
matter from fluid, MHD generators, pumps, accelerators and flow
meters. The results of these investigations are not applicable to the
flow of ionized gases. In an ionized gas where the density is low and/or
the magnetic field is very strong, the conductivity normal to the free
spiraling of electrons and ions about the magnetic lines of force before
suffering collisions; also, a current is induced in a direction normal
to both the electric and magnetic fields. The phenomena, well-known in
the literature, is called the Hall effect. The study of
magnetohydrodynamic flows with Hall currents has important engineering
applications in problems of magnetohydrodynamics generators and of Hall
accelerators as well as in flight magnetohydrodynamics
\cite{9,10,11,18,20,21,23,24}.

An understanding of the dynamics of fluids in porous media has practical
interest in such disparate fields as petroleum engineering and ground
water hydrology, with applications ranging from hydrocarbon migration in
reservoirs via packed-bed chemical reactors to agricultural drainage and
irrigation. In the widely used continuum approach to transport processes
in porous media, the differential equation governing the macroscopic
fluid motion is based on the experimentally established Darcy's law
\cite{2}, which accounts for the drag exerted by the porous medium.
Brinkman \cite{5} studied Darcy's law by adding to it a viscous term in
order to account for the vorticity diffusion caused by the boundary
resistance, whereas the combined influence of inertia and viscous
effects on the flow and heat transfer in the vicinity of an impermeable
plane surface was discussed by Vafai and Tien \cite{28} and Kaviany
\cite{13}. Kaviany's analysis was for steady-state flow, the oscillatory
motion was studied by Khodadadi \cite{14,15} and transient fluid motion
in a porous medium channel was discussed by Anderson and Holmedal
\cite{1}.

The governing equations that describe the flow of a Newtonian fluid is
the Navier--Stokes equations. These equations are non-linear partial
differential equations and known exact solutions are few in number.
Exact solutions are very important not only because they are solutions
of some fundamental flows but also because they serve as accuracy checks
for experimental, numerical and asymptotic methods. Since the equations
of motion of non-Newtonian fluids are more complicated and non-linear
than the Navier--Stokes equations, so the inverse methods described by
Nemenyi \cite{17} have become attractive. In these methods, solutions
are found by assuming certain physical or geometrical properties of the
flow field. Kaloni and Huschilt \cite{12}, Siddiqui and Kaloni
\cite{25}, Siddiqui \cite{26}, Benharbit and Siddiqui \cite{3},
Labropulu \cite{16} and Siddiqui {\it et~al} \cite{27} used this method
to study the flow problems of a second-grade fluid.

In this work, we discuss the effects of Hall currents on the steady flow
of an electrically conducting second-grade fluid in a porous medium
channel. For such a fluid equations are modeled for a grade of fluid two
and are solved by assuming certain form of the stream function. The
graphs are plotted explicitly in the functional form to see the
behaviour of the flow field.

The paper is arranged in the following fashion: In \S2, governing
equations and formulation of the problem are given. Section 3 consists
of two parts. First part is the generalization of Siddiqui's \cite{26}
work and the second part deals with some special flows called as
Riabounchinsky type flows and finally, concluding remarks are given in
\S4. Stream function, velocity components and the pressure fields are
derived in each case. Moreover, the streamlines are plotted in each case
to see the flow behaviour.

\section{Governing equations}

The constitutive equation of an incompressible fluid of second-grade is of
the form \cite{22}
\begin{equation}
\mathbf{T} = - p\mathbf{I} + \mu \mathbf{A}_{1} + \alpha_{1}
\mathbf{A}_{2} + \alpha_{2} \mathbf{A}_{1}^{2},
\end{equation}
where $\mathbf{T}$ is the Cauchy stress tensor, $-p\mathbf{I}$ denotes
the indeterminate spherical stress and $\mu, \alpha_{1}$ and
$\alpha_{2}$ are measurable material constants. They denote,
respectively, the viscosity, elasticity and cross-viscosity. These
material constants can be determined from viscometric flows for any real
fluid. $\mathbf{A}_{1}$ and $\mathbf{A}_{2}$ are Rivlin--Ericksen
tensors \cite{22} and they denote, respectively, the rate of strain and
acceleration. $\mathbf{A}_{1}$ and $\mathbf{A}_{2}$ are defined by
\begin{align}
\mathbf{A}_{1} &= (\hbox{grad}\mathbf{V}) + (\hbox{grad}\mathbf{V})^{\top},\\[.2pc]
\mathbf{A}_{2} &= \frac{\hbox{d}\mathbf{A}_{1}}{\hbox{d}t} +
\mathbf{A}_{1} (\hbox{grad} \mathbf{V}) +
(\hbox{grad}\mathbf{V})^{\top}\mathbf{A}_{1}.
\end{align}
Here $\mathbf{V}$ is the velocity, grad the gradient operator, $\top$
the transpose, and ${\rm d}/{\rm d}t$ the material time derivative.

The basic equations governing the motion of an incompressible fluid are
\begin{align}
&\mathbf{\nabla}\cdot \mathbf{V} = 0,\\[.2pc]
&\rho \frac{\hbox{d}\mathbf{V}}{\hbox{d}t} = \mathbf{J} \times \mathbf{B}
+\hbox{div}\,\mathbf{T}-\frac{\mu }{K} \mathbf{V},
\end{align}
\begin{align}
&\mathbf{\nabla }\cdot \mathbf{B} = 0,\text{ }\mathbf{\nabla }\times
\mathbf{B} = \mu_{m}\mathbf{J},\text{ }\mathbf{\nabla }\times \mathbf{E}
= 0,
\end{align}
where $\rho$ is the density, $\mathbf{J}$ the current density,
$\mathbf{B}$ the total magnetic field, $\mu_{m}$ the magnetic
permeability, $\mathbf{E}$ the total electric field current and $K$ the
permeability of the porous medium. Making reference to Cowling \cite{6},
when the strength of the magnetic field is very large, the generalized
Ohm's law is modified to include the Hall current so that
\begin{equation}
\mathbf{J} + \frac{\omega_{e}\tau_{e}}{\mathbf{B}_{0}} (\mathbf{J}
\times \mathbf{B}) = \sigma \left[ \mathbf{E} + \mathbf{V} \times
\mathbf{B} + \frac{1}{en_{e}}\mathbf{\nabla}p_{e}\right]
\end{equation}
in which $\omega_{e}$ is the cyclotron frequency, $\tau_{e}$ the
electron collision time, $\sigma$ the electrical conductivity, $e$
the electron charge and $p_{e}$ the electron pressure. The ion-slip and
thermoelectric effects are not included in (2.7). Further, it
is assumed that $\omega_{e}\tau_{e}\sim O1$ and $\omega_{i}\tau_{i} \ll 1$,
where $\omega_{i}$ and $\tau_{i}$ are the cyclotron frequency and
collision time for ions respectively.

Inserting (2.1) in (2.5) and making use of (2.2), (2.3), (2.6) and (2.7) we obtain the following vector equation
\begin{align}
&\hbox{grad} \left[\frac{1}{2}\rho |\mathbf{V}|^{2} + p - \alpha_{1}
\left( \mathbf{V\cdot \nabla}^{2} \mathbf{V} + \frac{1}{4}
|\mathbf{A}_{1}|^{2}\right) \right] + \rho [\mathbf{V}_{t} -
\mathbf{V\times} (\mathbf{\nabla} \times \mathbf{V})]\nonumber\\[.2pc]
&\quad\ = \mu \mathbf{\nabla}^{2} \mathbf{V} + \alpha_{1}
[\mathbf{\nabla}^{2} \mathbf{V}_{t} + \mathbf{\nabla}^{2}
(\mathbf{\nabla }\times \mathbf{V}) \times \mathbf{V}]\nonumber\\[.2pc]
&\qquad\ \ + (\alpha_{1} + \alpha_{2})\ \hbox{div}\ \mathbf{A}_{1}^{2} -
\frac{\sigma\mathbf{B}_{0}^{2} (1 + i\varphi)}{1 +\varphi^{2}} -
\frac{\mu}{K}\mathbf{V},
\end{align}
in which $\mathbf{\nabla}^{2}$ is the Laplacian operator,
$\mathbf{V}_{t} = \partial \mathbf{V}/\partial t$, and
$|\mathbf{A}_{1}|$ is the usual norm of matrix $\mathbf{A}_{1}$. If this
model is required to be compatible with thermodynamics, then the
material constants must meet the restrictions \cite{7,8}
\begin{equation}
\mu \geq 0,\text{ \ }\alpha_{1}\geq 0,\text{ \ }\alpha_{1} + \alpha_{2}
= 0.
\end{equation}

On the other hand, experimental results of tested fluids of second-grade
showed that $\alpha_{1} < 0$ and $\alpha_{1} + \alpha_{2}\neq 0$ which
contradicts the above conditions and imply that such fluids are
unstable. This controversy is discussed in detail in \cite{19}. However,
in this paper we will discuss both cases, $\alpha_{1}\geq 0$ and
$\alpha_{1} < 0$.

The velocity field for the problem under consideration is of the
following form
\begin{equation}
\mathbf{V} (x, y, t) = [u (x, y, t), \text{ }v (x, y, t),\text{ }0],
\end{equation}
where $u$ and $v$ are velocity components in the $x$ and $y$ directions,
respectively.

Inserting (2.10) in (2.4) and (2.8) and making use of the assumption
(2.9) we obtain the following equations
\begin{align}
&\frac{\partial u}{\partial x} + \frac{\partial v}{\partial y} = 0,\\[.2pc]
&\frac{\partial \widehat{p}}{\partial x} + \rho \left[ \frac{\partial
u}{\partial t} - v\omega \right] = \left( \mu +
\alpha_{1}\frac{\partial}{\partial t}\right)
\mathbf{\nabla}^{2}u - \alpha_{1}v\mathbf{\nabla}^{2}\omega\nonumber\\[.2pc]
&\hskip 8pc - \frac{\sigma \mathbf{B}_{0}^{2}(1 + i\varphi)}{1 +
\varphi^{2}} u - \frac{\mu}{K}u,
\end{align}
\begin{align}
\frac{\partial \widehat{p}}{\partial y} + \rho \left[ \frac{\partial
v}{\partial t} + u\omega \right] &= \left( \mu +
\alpha_{1}\frac{\partial}{\partial t}\right) \mathbf{\nabla}^{2}
v + \alpha_{1}u\mathbf{\nabla}^{2}\omega\nonumber\\[.2pc]
&\quad\ - \frac{\sigma \mathbf{B}_{0}^{2} (1 + i\varphi)}{1 +
\varphi^{2}} v - \frac{\mu}{K}v,
\end{align}
where $\varphi = \omega_{e}\tau_{e}$ is the Hall parameter and
\begin{align}
&\omega = \frac{\partial v}{\partial x} - \frac{\partial u}{\partial
y},\\[.2pc]
&\widehat{p} = p + \frac{1}{2}\rho (u^{2} + v^{2}) - \alpha_{1} \left[
u\mathbf{\nabla}^{2}u + v\mathbf{\nabla}^{2}v + \frac{1}{4}
|\mathbf{A}_{1}^{2}| \right],\\[.2pc]
&|\mathbf{A}_{1}^{2}| = 4\left( \frac{\partial u}{\partial x}
\right)^{2} + 4\left( \frac{\partial v}{\partial y}\right)^{2} + 2\left(
\frac{\partial u}{\partial y} + \frac{\partial v}{\partial
x}\right)^{2}.
\end{align}

\begin{rema}
{\rm On setting $\alpha _{1}=0, K\rightarrow \infty$ and neglecting Hall
effects in (2.12) and (2.13) we recover the equations for viscous fluid,
on taking $K\rightarrow \infty$ and neglecting Hall effects we obtain
the case \cite{26} and on taking $\alpha_{1} = 0$ and $\varphi = 0$ in
Hall effects we obtain the Brinkman model for porous medium.}
\end{rema}

Equations (2.11)--(2.13) are three partial differential equations for
three unknown functions $u, v$ and $\widehat{p}$ of the variables
$(x,y)$. Once the velocity field is determined, the pressure field
(2.15) can be calculated by integrating (2.12) and (2.13). Note that the
equation for the component $w$ is identically zero.

Eliminating pressure in (2.12) and (2.13), by applying the integrability
condition $\partial^{2}\widehat{p}/\partial x\partial y =
\partial^{2}\widehat{p}/\partial y\partial x$, we get the compatibility
equation
\begin{align}
\rho \left[ \frac{\partial \omega}{\partial t} + \left(
u\frac{\partial}{\partial x} + v\frac{\partial}{\partial y}\right)
\omega \right] &= \left( \mu + \alpha_{1}\frac{\partial}{\partial
t}\right) \mathbf{\nabla}^{2}\omega - \frac{\sigma \mathbf{B}_{0}^{2} (1
+ i\varphi)}{1 + \varphi^{2}}\omega\nonumber\\[.2pc]
&\quad\ -\frac{\mu}{K}\omega + \alpha_{1}\left[ \left(
u\frac{\partial}{\partial x} + v\frac{\partial}{\partial y}\right)
\mathbf{\nabla}^{2}\omega \right].
\end{align}
Let us consider the Stokes stream function:
\begin{equation}
u = \frac{\partial \psi}{\partial y},\text{ \ } v = -\frac{\partial
\psi}{\partial x},
\end{equation}
where $\psi (x,y)$ is the stream function. We see that the continuity
equation (2.11) is satisfied identically and (2.18) in (2.17) yields the
following equation
\begin{align}
\rho \left[ \frac{\partial}{\partial t}\mathbf{\nabla}^{2}\psi - \{\psi,
\mathbf{\nabla}^{2}\psi\} \right] &= \left( \mu +
\alpha_{1}\frac{\partial}{\partial t}\right) \mathbf{\nabla}^{4}\psi -
\alpha_{1} [\{\psi, \mathbf{\nabla}^{4}\psi\}]\nonumber\\[.2pc]
&\quad\ - \frac{\sigma \mathbf{B}_{0}^{2} (1 + i\varphi)}{1 +
\varphi^{2}}\mathbf{\nabla}^{2}\psi - \frac{\mu}{K}\mathbf{\nabla}^{2}
\psi,
\end{align}
in which
\begin{equation}
\mathbf{\nabla}^{4} = \mathbf{\nabla}^{2}\cdot
\mathbf{\nabla}^{2},\text{ \ }\omega =- \mathbf{\nabla}^{2}\psi
\end{equation}
and
\begin{equation}
\{\psi,\mathbf{\nabla}^{2}\psi\} = \frac{\partial \psi}{\partial
x}\frac{\partial (\mathbf{\nabla}^{2}\psi)}{\partial y} - \frac{\partial
\psi}{\partial y}\frac{\partial (\mathbf{\nabla}^{2}\psi)}{\partial x}
\end{equation}
is the Poisson bracket.

\begin{rema}
{\rm The equation (2.19), for $K\rightarrow \infty$ and in the absence
of Hall effects reduces to~\cite{26}.}
\end{rema}

\section{Solutions of some special types}

\subsection{\it Solution of the type $\protect\psi (x,y) = \protect\xi
(x) + \protect\eta (y)$}

We consider the plane steady flow and examine the solution of (2.19) of
the form
\renewcommand{\theequation}{\thesubsection.\arabic{equation}}
\setcounter{equation}{0}
\begin{equation}
\psi (x,y) = \xi (x) + \eta (y),
\end{equation}
where $\xi$ and $\eta$ are arbitrary functions of the variables $x$ and
$y$ respectively. Substituting (3.1.1) in (2.19) we obtain the following
equation
\begin{align}
\rho [\eta^{\prime} (y) \xi^{\prime \prime \prime} (x) - \xi^{\prime}
(x) \eta^{\prime \prime \prime} (y)] &= \mu [\xi^{\rm IV} (x) +
\eta^{\rm IV} (y)]\nonumber\\[.2pc]
&\quad\ + \alpha_{1} [\eta^{\prime} (y) \xi^{\rm V} (x) - \xi^{\prime}
(x) \eta^{\rm V} (y)]\nonumber\\[.2pc]
&\quad\ + \left(\frac{\sigma \mathbf{B}_{0}^{2} (1 + i\varphi)}{1 +
\varphi^{2}} + \frac{\mu}{K}\right) [\xi^{\prime \prime} (x) +
\eta^{\prime \prime} (y)],
\end{align}
in which IV and V in the superscript indicates the fourth and fifth
derivatives.

We see that (3.1.2) is highly non-linear and its solution in the present
form is not easy to obtain. In order to find its solution we assume the
following
\begin{align}
\xi (x) &= Ax + B{\rm e}^{ax},\\[.2pc]
\eta (y) &= Cy + D{\rm e}^{by}
\end{align}
and obtain the following equation
\begin{align}
&\rho [-Ab^{3} D{\rm e}^{by} + a^{3} BC{\rm e}^{ax} - ab^{3} BD{\rm e}^{ax +
by} + a^{3} bBD{\rm e}^{ax + by}]\nonumber\\[.2pc]
&\quad\ = \mu [a^{4} B{\rm e}^{ax} + b^{4} D{\rm e}^{by}]
+ \alpha_{1} \left[ \begin{array}{@{}c@{}}
- Ab^{5} D{\rm e}^{by} + a^{5} BC{\rm e}^{ax}\\[.3pc]
- ab^{5} BD{\rm e}^{ax + by} + a^{5}bBD{\rm e}^{ax + by}
\end{array}\right]\nonumber\\[.2pc]
&\qquad\ - \left( \frac{\sigma \mathbf{B}_{0}^{2} (1 + i\varphi)}{1 +
\varphi^{2}} + \frac{\mu}{K}\right) [a^{2} B{\rm e}^{ax} + b^{2}D{\rm
e}^{by}],
\end{align}
where $A, B, C, D, a$ and $b$ are arbitrary constants.

The following three equations are obtained from (3.1.5)
\begin{align}
&\rho aC = \mu a^{2} + \alpha_{1}a^{3}C - H - \frac{\mu}{K},\\[.2pc]
&-\rho bA = \mu b^{2} - \alpha_{1}b^{3}A - H - \frac{\mu}{K},\\[.2pc]
&(b^{2} - a^{2}) [\rho - \alpha_{1} (a^{2} + b^{2})] = 0.
\end{align}

From (3.1.6) and (3.1.7) we easily obtain the values of $A$ and $C$, i.e.
\begin{align}
A &= - \frac{1}{\rho -\alpha_{1}b^{2}} \left[ \mu b - \frac{H}{b} -
\frac{\mu}{Kb} \right],\\[.3pc]
C &= \frac{1}{\rho -\alpha_{1}a^{2}} \left[ \mu a - \frac{H}{a} -
\frac{\mu}{Ka} \right],
\end{align}
where
\begin{equation*}
H = \frac{\sigma \mathbf{B}_{0}^{2} (1 + i\varphi)}{1 + \varphi^{2}}
\end{equation*}
and (3.1.8) is satisfied if either
\begin{equation}
b^{2} - a^{2} = 0
\end{equation}
or
\begin{equation}
\rho = \alpha_{1} (a^{2} + b^{2}).
\end{equation}

We have three different cases which we discuss separately as follows:

\begin{case}{\rm
$b = a, \rho \neq \alpha_{1} (a^{2} + b^{2})$

The stream function given by (3.1.1), after using (3.1.3), (3.1.4),
(3.1.9), and (3.1.10) becomes
\begin{equation}
\psi (x,y) = \frac{(y - x)}{\rho - \alpha_{1}a^{2}} \left[ \mu a -
\frac{H}{a} - \frac{\mu}{Ka}\right] + B{\rm e}^{ax} + D{\rm e}^{ay}
\end{equation}
and from (2.18) the velocity components take the following form
\begin{align}
u &= \frac{1}{\rho -\alpha_{1}a^{2}}\left[ \mu a - \frac{H}{a} -
\frac{\mu}{Ka} \right] + Da{\rm e}^{ay},\\[.2pc]
v &= \frac{1}{\rho -\alpha_{1}a^{2}}\left[ \mu a - \frac{H}{a} -
\frac{\mu}{Ka} \right] - Ba{\rm e}^{ax}.
\end{align}
In order to find the pressure field (2.15) we substitute the velocity
components (3.1.14) and (3.1.15) in (2.12) and (2.13) and then integrate
the resulting equations to obtain
\begin{align}
p &= p_{0} - \rho \overline{a}^{2} - \mu Ba^{3}y{\rm e}^{ax} + (\rho -
\alpha_{1}a^{2}) [a^{2}By{\rm e}^{ax} + a^{2}DB{\rm e}^{a(x +
y)}]\nonumber\\[.2pc]
&\quad\ + \alpha_{1} [B^{2}a^{4}{\rm e}^{2ax} + D^{2}a^{4}{\rm e}^{2ay}
- DBa^{4}{\rm e}^{a (x + y)}],
\end{align}
where
\begin{equation*}
\overline{a} = \frac{1}{\rho - \alpha_{1}a^{2}} \left[ \mu a -
\frac{H}{a} - \frac{\mu}{Ka}\right],
\end{equation*}
whereas the streamline for $\psi = \Omega_{1}$ is given by the following
functional form
\begin{equation}
y = \frac{- B{\rm e}^{ax} + x\varepsilon + \Omega_{1}}{\varepsilon} -
\frac{1}{a}\ \ \hbox{Product Log}\ \ \left[ \frac{Da}{\varepsilon} {\rm
e}^{a(-B{\rm e}^{ax} + x\varepsilon + \Omega_{1})/\varepsilon}\right],
\end{equation}
where
\begin{equation*}
\varepsilon = \frac{1}{1-\Lambda a^{2}}\left[ \nu a - \frac{\chi}{a} -
\frac{\nu}{Ka}\right],
\end{equation*}
with $\nu = \mu/\rho$ as the kinematic viscosity, $\Lambda =
\alpha_{1}/\rho$ is the second-grade parameter $\chi = N (1 +
i\varphi)/(1 + \varphi^{2})$, $N = \sigma \mathbf{B}_{0}^{2}/\rho$ is
the MHD parameter and $\varphi = \omega_{e}\tau_{e}$ is the Hall
parameter.

Streamlines are shown in figure~1 for $B = D = a = 1$, $\mu /\rho =
0.5$, $\alpha_{1}/\rho = 0.1$, $K = 0.1$, $N = 0$, $\omega_{e}\tau_{e} =
0.1$ and $\psi = 15, 20, 25, 30, 40$.}
\end{case}

\begin{fig}
\hskip 4pc{\epsfxsize=8.7cm\epsfbox{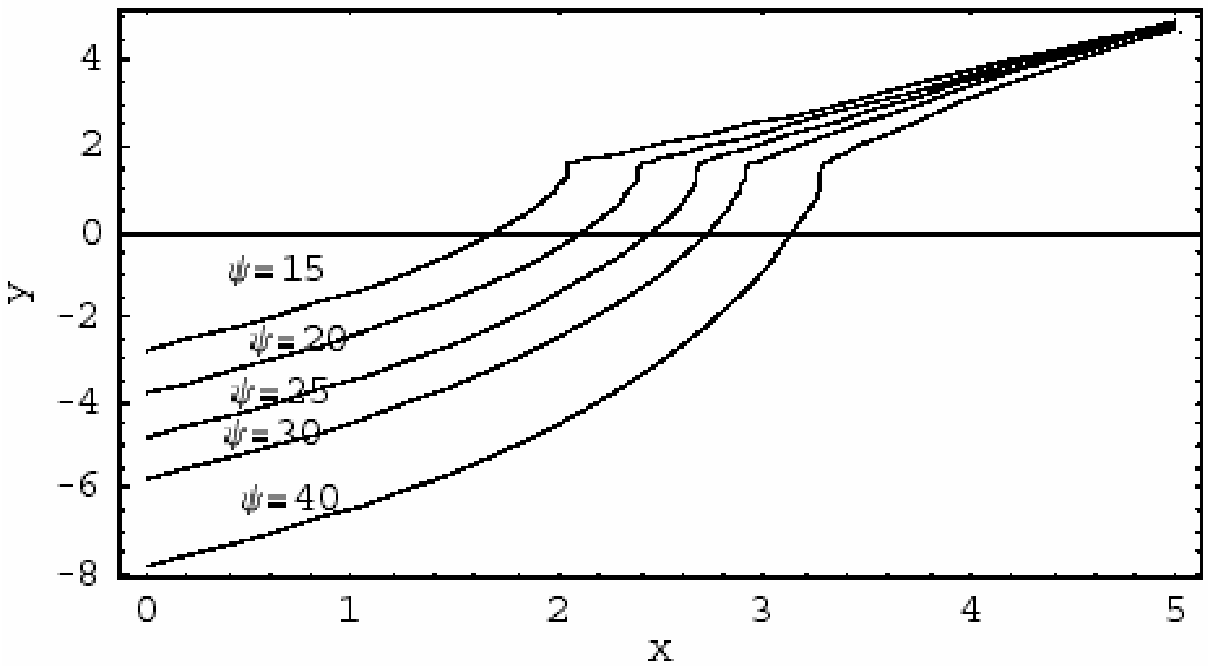}}\vspace{-.7pc}
\caption{Streamline flow pattern for $\protect\psi (x, y) = \frac{(y -
x)}{\rho - \alpha_{1}a^{2}} \left[\mu a - \frac{H}{a} - \frac{\mu
}{Ka}\right] + B{\rm e}^{ax}$\break}
\hskip 4pc{\xxxx $ + D{\rm e}^{ay}$.}\vspace{.5pc}
\end{fig}

\begin{case}{\rm
$b = -a, \rho \neq \alpha_{1} (a^{2} + b^{2})$

The expressions for $\psi, u, v$, and $p$ are
\begin{align}
&\psi (x, y) = \frac{1}{\rho -\alpha_{1}a^{2}}\left[ \mu a - \frac{H}{a}
- \frac{\mu}{Ka}\right] (y + x) + B{\rm e}^{ax} + D{\rm e}^{-ay},\\[.2pc]
&u = \frac{1}{\rho -\alpha_{1}a^{2}} \left[ \mu a - \frac{H}{a} -
\frac{\mu}{Ka} \right] - Da{\rm e}^{-ay},\\[.2pc]
&v = -\frac{1}{\rho -\alpha_{1}a^{2}}\left[ \mu a - \frac{H}{a} -
\frac{\mu}{Ka} \right] - Ba{\rm e}^{ax},\\[.2pc]
&p = p_{0} - \rho \overline{a}^{2} - \mu Ba^{3} y{\rm e}^{ax} + (\rho -
\alpha_{1} a^{2}) [a^{2}By{\rm e}^{ax} + a^{2}DB{\rm e}^{a (x - y)}]\nonumber\\[.2pc]
&\hskip 1.5pc + \alpha_{1} [B^{2}a^{4}{\rm e}^{2ax} + D^{2}a^{4}{\rm e}^{-2ay}
- DBa^{4}{\rm e}^{a(x - y)}]
\end{align}
and the functional form of streamline for $\psi = \Omega_{2}$ is given
by
\begin{equation}
y = \frac{-B{\rm e}^{ax} - x\varepsilon + \Omega_{2}}{\varepsilon} +
\frac{1}{a}\ \ \hbox{Product Log}\ \ \left[ -\frac{Da}{\varepsilon}{\rm
e}^{a (B{\rm e}^{ax} + x\varepsilon - \Omega_{2})/\varepsilon}\right].
\end{equation}
Streamlines are drawn in figure~2 for $B = D = a = 1$, $\mu/\rho = 0.5$,
$\alpha_{1}/\rho = 0.1$, $K = 0.1$, $N = 0.5$, $\omega_{e}\tau_{e} = 1$
and $\psi = 15, 20, 25, 30, 40$.}
\end{case}

\begin{fig}
\hskip 4pc{\epsfxsize=8.7cm\epsfbox{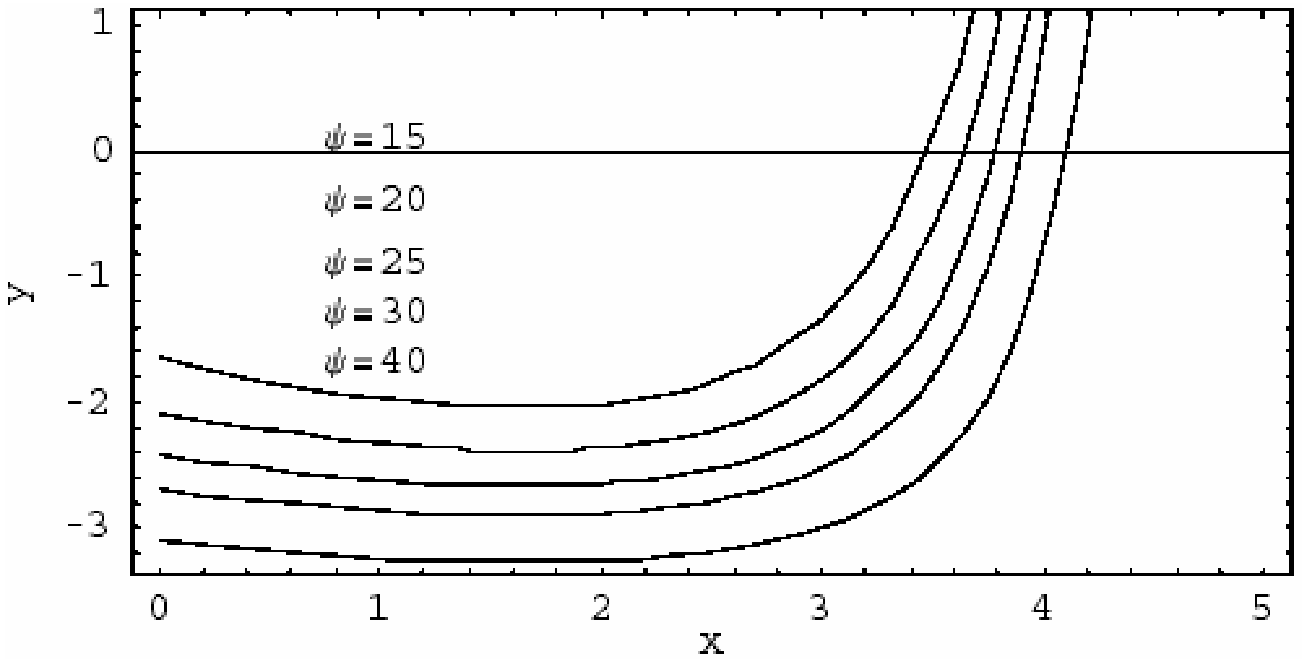}}\vspace{-.7pc}
\caption{Streamline flow pattern for $\psi (x, y) = \frac{1}{\rho -
\alpha_{1}a^{2}}\left[ \mu a - \frac{H}{a} - \frac{\mu}{Ka} \right] (y +
x)$\break}
\hskip 4pc {\xxxx $+ B{\rm e}^{ax} + D{\rm e}^{-ay}$.}\vspace{.5pc}
\end{fig}

\begin{case}{\rm
$b^{2} - a^{2}\neq 0$

We must have
\begin{equation*}
\rho = \alpha_{1} (a^{2} + b^{2})
\end{equation*}
and the expressions for $\psi, u, v$, and $p$ are of the following form
\begin{align}
&\psi (x, y) = -\frac{x}{\rho - \alpha_{1}b^{2}} \left[\mu b -
\frac{H}{b} - \frac{\mu}{Kb}\right]\nonumber\\[.2pc]
&\hskip 3.6pc + \frac{y}{\rho - \alpha_{1}a^{2}} \left[\mu a - \frac{H}{a} -
\frac{\mu}{Ka} \right] + B{\rm e}^{ax} + D{\rm e}^{by},\\[.2pc]
&u = \frac{1}{\rho - \alpha_{1}a^{2}} \left[\mu a - \frac{H}{a} -
\frac{\mu}{Ka} \right] + Db{\rm e}^{by},\\[.2pc]
&v = \frac{1}{\rho - \alpha _{1}b^{2}} \left[\mu b - \frac{H}{b} -
\frac{\mu}{Kb} \right] - Ba{\rm e}^{ax},\\[.2pc]
&p = p_{0} - \frac{1}{2}\rho [\overline{a}^{2} + \overline{b}^{2}] - \mu
Ba^{3}y{\rm e}^{ax}\nonumber\\[.2pc]
&\hskip 1.4pc + (\rho - \alpha_{1}a^{2}) [a^{2} By{\rm e}^{ax} +
a^{2}DB{\rm e}^{ax + by}]\nonumber\\[.2pc]
&\hskip 1.4pc + \alpha_{1} [B^{2}a^{4}{\rm e}^{2ax} + D^{2}a^{4}{\rm
e}^{2by} - DBa^{4}b^{2}{\rm e}^{ax + by}],
\end{align}
where
\begin{equation*}
\overline{a} = \frac{1}{\rho - \alpha_{1}a^{2}}\left[ \mu a -
\frac{H}{a} - \frac{\mu}{Ka}\right],\quad\overline{b} = \frac{1}{\rho
- \alpha_{1}b^{2}} \left[\mu b - \frac{H}{b} - \frac{\mu}{Kb}\right],
\end{equation*}
whereas the functional form in this case for $\psi = \Omega_{3}$ is
\begin{equation}
y = -\frac{B{\rm e}^{ax} - x\delta - \Omega_{3}}{\varepsilon} -
\frac{1}{b} \ \ \hbox{Product Log}\ \ \left[ \frac{Db}{\varepsilon} {\rm
e}^{-b (B{\rm e}^{ax} - x\delta - \Omega_{3})/\varepsilon}\right],
\end{equation}
where
\begin{equation*}
\delta = \frac{1}{1 - \Lambda b^{2}} \left[\nu b - \frac{\chi}{b} -
\frac{\nu}{Kb} \right].
\end{equation*}
Streamlines for $B = D = a = 1$, $\mu/\rho = 0.5$, $b = -0.5$, $K = 2,
2.1, \alpha_{1}/\rho = 0.5$, $-0.7, N = 0.5$, $\omega_{e}\tau_{e} = 1$
and $\psi = 15, 20, 25, 30, 40$ are depicted in figures~3, 4, whereas
figure~5 is given for $\alpha_{1}/\rho =- 0.5$.

\begin{fig}[b]
\hskip 4pc{\epsfxsize=8.7cm\epsfbox{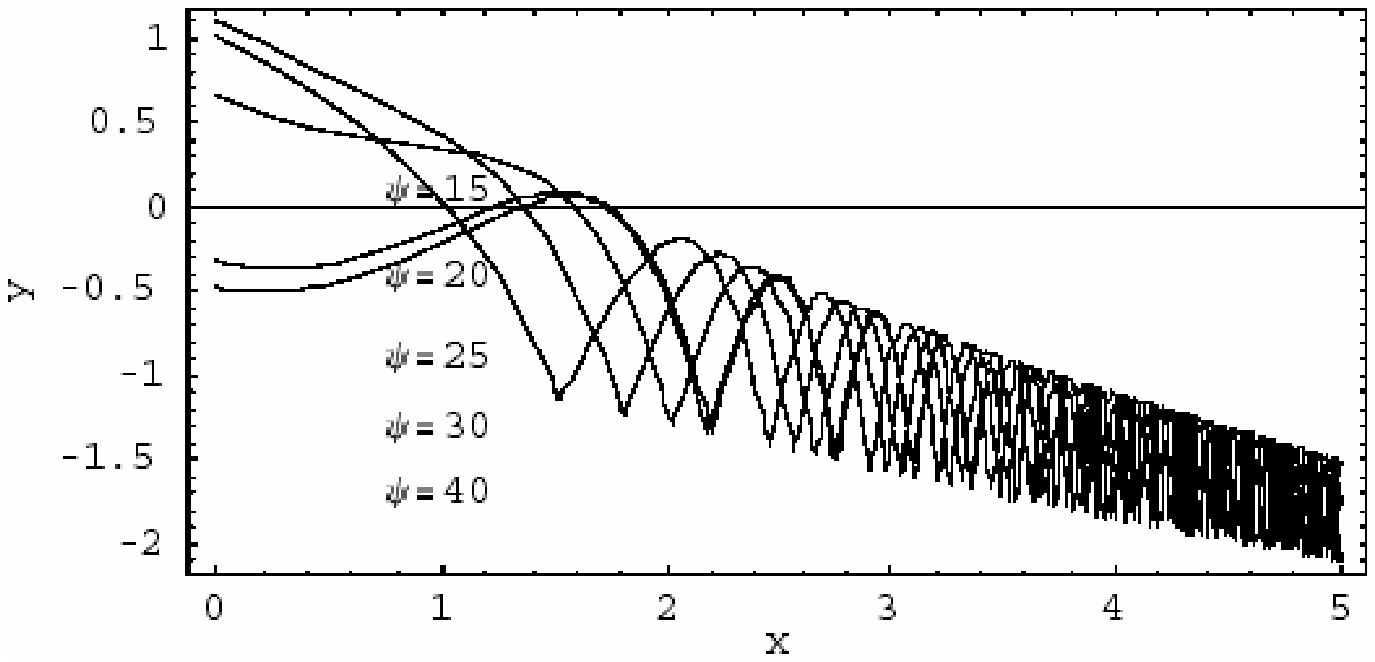}}\vspace{-.5pc}
\caption{$\psi (x, y) = -x\overline{a} + y\overline{b} + B{\rm e}^{ax} +
D{\rm e}^{by}$.}\vspace{.5pc}
\end{fig}

\begin{fig}
\hskip 4pc{\epsfxsize=8.7cm\epsfbox{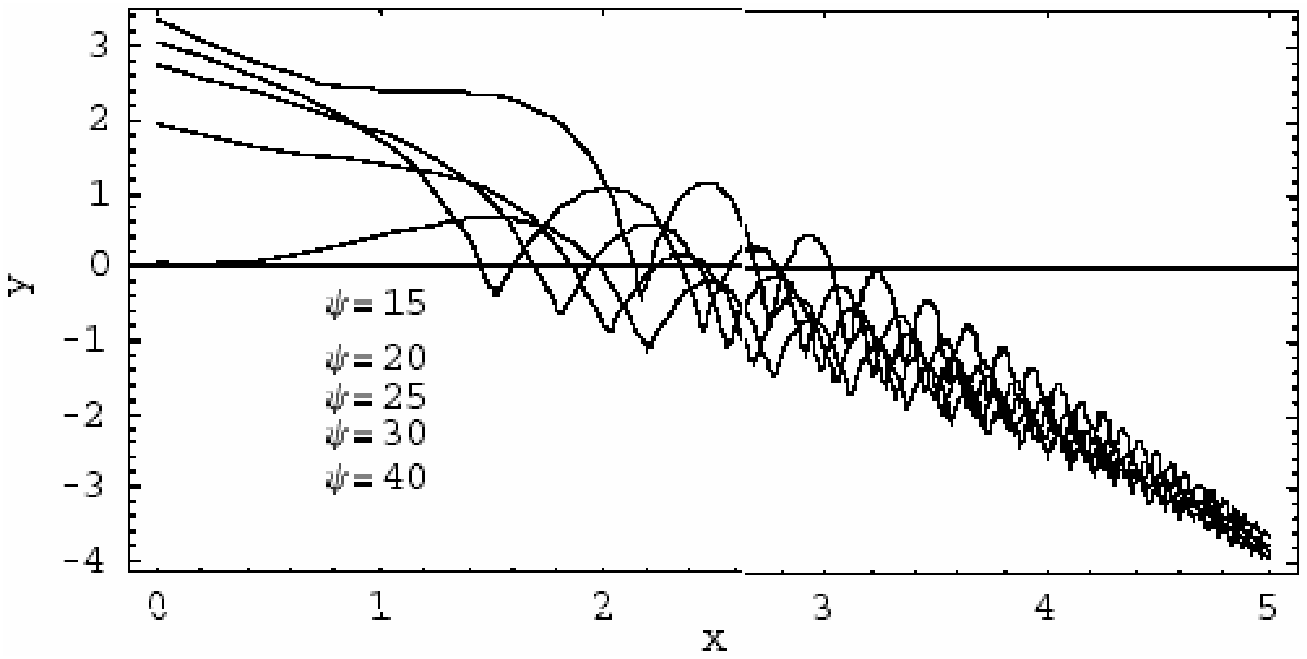}}\vspace{-.5pc}
\caption{Streamline flow pattern for $\psi (x, y) = - x\overline{a} +
y\overline{b} + B{\rm e}^{ax} + D{\rm e}^{by}$.}\vspace{.5pc}
\end{fig}

\begin{fig}
\hskip 4pc{\epsfxsize=8.7cm\epsfbox{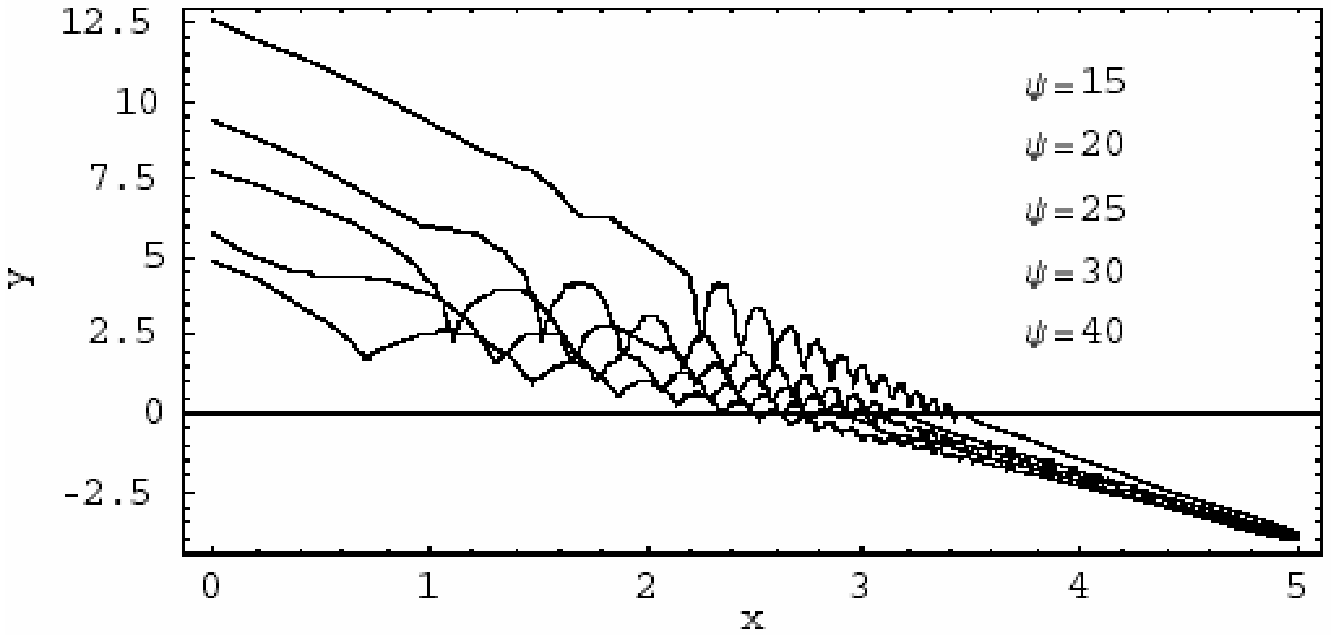}}\vspace{-.5pc}
\caption{Streamline flow pattern for negative second-grade parameter.}
\end{fig}

The alternate forms of (3.1.23) may be written as
\begin{align}
\psi (x, y) &= -\frac{x}{\rho -\alpha_{1}b^{2}} \left[ \mu b -
\frac{H}{b} - \frac{\mu}{Kb} \right]\nonumber\\[.2pc]
&\quad\ + \frac{y}{\alpha_{1}b^{2}} \left[ \mu a -\frac{H}{a} -
\frac{\mu}{Ka}\right] + B{\rm e}^{ax} + D{\rm e}^{by},\\[.2pc]
\psi (x, y) &= -\frac{x}{\alpha_{1} a^{2}} \left[\mu b - \frac{H}{b} -
\frac{\mu}{Kb} \right]\nonumber\\[.2pc]
&\quad\ + \frac{y}{\rho - \alpha_{1}a^{2}} \left[\mu a - \frac{H}{a} -
\frac{\mu}{Ka} \right] + B{\rm e}^{ax} + D{\rm e}^{by}.
\end{align}}
\end{case}

\begin{rema}{\rm
The solution (3.1.16) with $H = 0$, $K\rightarrow \infty$ and
$\alpha_{1} = 0$ gives the Berker's solution \cite{4} and the Siddiqui's
solutions \cite{26} can readily be recovered as a special case by taking
$H = 0$ and $K\rightarrow \infty$.}
\end{rema}

We now consider {\it Riabounchinsky type flows} in order to solve (2.15).

\subsection{\it Solution of the type $\protect\psi (x, y) = y\xi (x)$}

In order to obtain another class of solution of (2.19) we substitute
\setcounter{equation}{0}
\begin{equation}
\psi (x, y) = y\xi (x)
\end{equation}
into (2.19) and get the following equation
\begin{equation}
\rho [\xi \xi^{\prime \prime \prime} - \xi^{\prime} \xi^{\prime\prime}]
= \mu \xi^{\rm IV} + \alpha_{1} [\xi \xi^{\rm V} - \xi^{\prime} \xi^{\rm
IV}] - \left( H + \frac{\mu}{K}\right) \xi^{\prime\prime},
\end{equation}
where $\xi (x)$ is an arbitrary function of $x$, primes denote the
derivative with respect to $x$.

Integrating (3.2.2) once and equating the constant of integration equal
to zero we obtain
\begin{equation}
\mu \xi^{\prime\prime\prime} + \rho [(\xi^{\prime 2} - \xi\xi^{\prime
\prime})] + \alpha_{1} [(\xi \xi^{\rm IV} - 2\xi^{\prime} \xi^{\prime
\prime \prime} + \xi^{\prime \prime 2})] - \left( H +
\frac{\mu}{K}\right) \xi^{\prime} = 0.
\end{equation}

For the solution of the above equation we write
\begin{equation}
\xi (x) = \delta (1 + \lambda {\rm e}^{\sigma x})
\end{equation}
in which $\delta, \sigma$ and $\lambda$ are arbitrary real constants.
Making use of (3.2.4) into (3.2.3) we have
\begin{equation}
\delta = \frac{1}{\rho - \alpha _{1}\sigma^{2}}\left[ \mu \sigma -
\frac{1}{\sigma} \left( H + \frac{\mu}{K}\right) \right]
\end{equation}
and thus from (3.2.1)
\begin{equation}
\psi (x, y) = \frac{1}{\rho - \alpha_{1}\sigma^{2}}\left[ \mu \sigma -
\frac{1}{\sigma} \left( H + \frac{\mu}{K}\right) \right] y (1 + \lambda
{\rm e}^{\sigma x}).
\end{equation}
The velocity components (2.18) and the pressure field (2.15) become
\begin{align}
u &= \frac{1}{\rho - \alpha_{1}\sigma^{2}} \left[ \mu \sigma -
\frac{1}{\sigma} \left( H + \frac{\mu}{K}\right) \right] (1 + \lambda
{\rm e}^{\sigma x}),\\[.2pc]
v &= \frac{-y}{\rho - \alpha_{1}\sigma^{2}} \left[ \mu \sigma -
\frac{1}{\sigma} \left( H + \frac{\mu}{K}\right) \right] \sigma \lambda
{\rm e}^{\sigma x}.\\[.2pc]
p &= p_{1} + \mu \sigma \overline{a} \left( 1 - \frac{\sigma^{2}
y^{2}}{2}\right) \lambda {\rm e}^{\sigma x} - \frac{1}{2}\rho
\left[ \begin{array}{@{}c@{}}
\overline{a}^{2} (1 - \lambda^{2}{\rm e}^{2\sigma x})
\end{array} \right]\nonumber\\[.2pc]
&\quad\ + \alpha_{1} \left[ \overline{a}^{2} \sigma^{2}\lambda {\rm
e}^{\sigma x} + \overline{a}^{2} \sigma^{2}\lambda^{2}\left( 3 +
\frac{\sigma^{2}y^{2}}{2}\right) {\rm e}^{2\sigma x} \right],
\end{align}
where $p_{1}$ is the reference pressure and
\begin{equation*}
\overline{a} = \frac{1}{\rho - \alpha_{1}\sigma^{2}}\left[ \mu \sigma -
\frac{1}{\sigma}\left( H + \frac{\mu}{K}\right) \right].
\end{equation*}
The streamline flow for $\psi = \Omega_{4}$ is given by the functional
form
\begin{equation}
y = \frac{\Omega_{4}}{(1 + \lambda {\rm e}^{\sigma x}) \varepsilon},
\end{equation}
where
\begin{equation*}
\varepsilon = \frac{1}{1 - \Lambda \sigma^{2}}\left[ \nu \sigma -
\frac{1}{\sigma} \left(\chi + \frac{\nu}{K} \right) \right].
\end{equation*}
Figure~6 shows the streamlines for $\sigma = \lambda = 1$, $\mu/\rho =
0.5$, $\alpha_{1}/\rho = 0.1$, $K = 15, N = 0, \psi = 15, 20, 25, 30,
40$.

\begin{fig}
\hskip 4pc{\epsfxsize=8.7cm\epsfbox{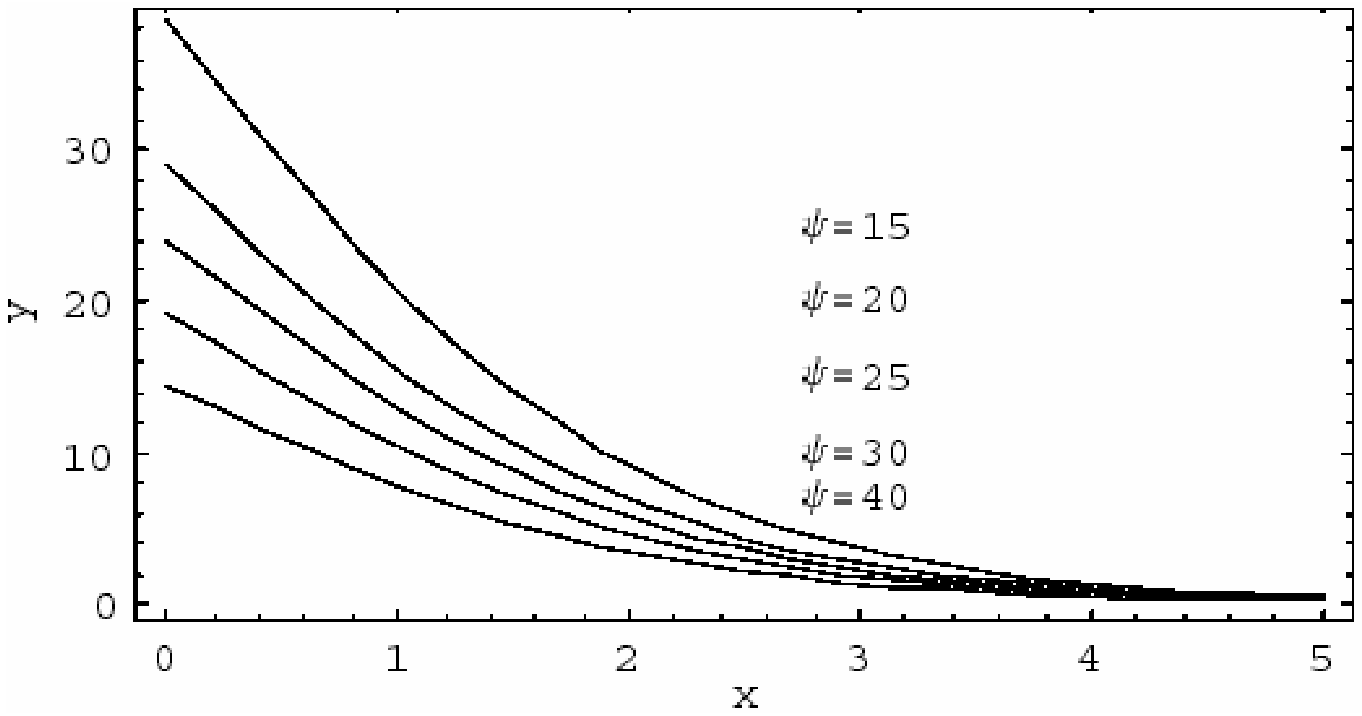}}\vspace{-.7pc}
\caption{Streamline flow pattern for $\protect\psi (x, y) =
\frac{1}{\rho - \protect \alpha_{1}\sigma^{2}} \left[ \protect\mu
\protect\sigma - \frac{1}{\sigma}\left( H + \frac{\mu}{K} \right)
\right] y(1 + $}
\hskip 4pc{\xxxx $\lambda {\rm e}^{\sigma x})$.}\vspace{.4pc}
\end{fig}

\subsection{\it Solutions of the type $\protect \psi(x, y) = y
\protect\xi(x) + \protect\eta (x)$}

Inserting
\setcounter{equation}{0}
\begin{equation}
\psi (x, y) = y\xi (x) + \eta (x)
\end{equation}
in (2.19) we obtain the following equation
\begin{align}
&- \rho [y (\xi^{\prime} \xi^{\prime \prime} - \xi \xi^{\prime \prime
\prime}) + (\eta^{\prime }\xi ^{\prime \prime} - \xi \eta^{\prime \prime
\prime})]\nonumber\\[.2pc]
&\quad\ = \mu (y\xi^{\rm IV} + \eta^{\rm IV}) - \alpha_{1} [y
(\xi^{\prime}\xi^{\rm IV} - \xi\xi^{\rm V}) +
(\eta^{\prime} \xi^{\rm IV} - \xi \eta^{\rm V})]\nonumber\\[.2pc]
&\qquad\ - \left( H + \frac{\mu}{K}\right) (y\xi^{\prime \prime} +
\eta^{\prime \prime}).
\end{align}
From the above equation, we have
\begin{equation}
\rho [\xi^{\prime}\xi^{\prime \prime} - \xi \xi^{\prime \prime \prime}]
+ \mu \xi^{\rm IV} - \alpha_{1} [\xi^{\prime}\xi^{\rm IV} - \xi \xi^{\rm
V}] - \left( H + \frac{\mu}{K}\right) \xi^{\prime \prime} = 0
\end{equation}
and
\begin{equation}
\rho [\eta^{\prime}\xi^{\prime \prime} - \xi \eta^{\prime \prime
\prime}] + \mu \eta^{\rm IV} - \alpha_{1} [\eta^{\prime} \xi^{\rm IV} -
\xi \eta^{\rm V}] - \left( H + \frac{\mu}{K}\right) \eta^{\prime \prime}
= 0,
\end{equation}
where $\xi (x)$ and $\eta (x)$ are arbitrary functions of its arguments.
Integrating (3.3.3) and (3.3.4) and then taking the constants of
integration equal to zero we have
\begin{align}
&\mu \xi^{\prime \prime \prime} + \rho [\xi^{\prime 2} - \xi \xi^{\prime
\prime}] - \alpha_{1} [(-\xi \xi^{\rm IV} + 2\xi^{\prime} \xi^{\prime
\prime \prime} - \xi^{\prime \prime 2})] - \left( H +
\frac{\mu}{K}\right) \xi^{\prime} = 0,\\[.2pc]
&\mu \eta^{\prime \prime \prime} + \rho [\eta^{\prime}\xi^{\prime} - \xi
\eta^{\prime \prime}] - \alpha_{1}
\left[ \begin{array}{@{}c@{}}
\xi^{\prime}\eta^{\prime \prime \prime} - \xi \eta^{\rm IV}\\[.2pc]
+ \eta^{\prime} \xi^{\prime \prime \prime} - \eta^{\prime \prime}
\xi^{\prime \prime}
\end{array} \right] - \left( H + \frac{\mu}{K}\right) \eta^{\prime} = 0.
\end{align}

We note that (3.3.5) is similar to (3.2.3). Its solution is given in
(3.2.4). Substituting (3.2.4) into (3.3.6) we have
\begin{align}
&\alpha_{1} A (1 + \lambda {\rm e}^{\sigma x}) \eta^{\rm IV} + (\mu
- \alpha_{1} A \lambda \sigma {\rm e}^{\sigma x}) \eta^{\prime \prime
\prime} + A [(\alpha_{1}\sigma^{2} - \rho) \lambda {\rm e}^{\sigma x} -
\rho] \eta^{\prime \prime}\nonumber\\[.2pc]
&\quad\ = -\left[ (\rho - \alpha_{1}\sigma^{2}) A\lambda \sigma {\rm
e}^{\sigma x} - \left( H + \frac{\mu}{K}\right) \right] \eta^{\prime},
\end{align}
where
\begin{equation*}
A = \frac{1}{\rho -\alpha_{1}\sigma^{2}} \left[ \mu\sigma -
\frac{1}{\sigma} \left( H + \frac{\mu}{K}\right) \right].
\end{equation*}

We note that it is not easy to obtain the general solution of (3.3.7).
In order to find its solution we consider the following special cases:

\setcounter{case}{0}
\begin{case}{\rm $\alpha _{1}\neq 0, \sigma =1, \lambda =0$

Equation (3.3.7) reduces to
\begin{equation}
\alpha_{1} A_{1} \eta^{\rm IV} + \mu\eta^{\prime \prime \prime} - \rho
A_{1} \eta^{\prime \prime} - \left( H + \frac{\mu}{K}\right)
\eta^{\prime} = 0.
\end{equation}

We see that (3.3.8) is of fifth order and in order to solve it we reduce
its order by putting $\eta^{\prime} = \overline{A} (x)$ such that
(3.3.8) becomes
\begin{equation}
\alpha_{1} A_{1} \overline{A}^{\prime \prime \prime} + \mu
\overline{A}^{\prime \prime} - \rho A_{1}\overline{A}^{\prime} - \left(
H + \frac{\mu}{K} \right) \overline{A} = 0.
\end{equation}
On substituting $\overline{A} (x) = \widehat{P} (x) {\rm e}^{x}$,
(3.3.9) takes the form
\begin{equation}
\alpha_{1} (3\widehat{P}^{\prime} + 3\widehat{P}^{\prime \prime}+
\widehat{P}^{\prime \prime \prime}) {\rm e}^{x} + \frac{\mu}{A_{1}}
(2\widehat{P}^{\prime} + \widehat{P}^{\prime \prime}) {\rm e}^{x} - \rho
\widehat{P}^{\prime} {\rm e}^{x} = 0.
\end{equation}
Finally, $\widehat{P}^{\prime}(x) = R(x)$ converts (3.3.10) into a
second-order differential equation
\begin{equation}
\alpha_{1} R^{\prime \prime} + (\mu /A_{1} + 3\alpha_{1}) R^{\prime} +
(3\alpha_{1} - \rho + 2\mu/A_{1}) R = 0.
\end{equation}

The solution of the above equation is
\begin{equation}
R (x) = A_{3}\exp \left( \frac{- c - \sqrt{c^{2} - 4d}}{2}\right) x +
A_{4}\exp \left( \frac{- c + \sqrt{c^{2} - 4d}}{2}\right) x,
\end{equation}
where $A_{3}$ and $A_{4}$ are arbitrary constants and
\begin{align*}
&c = \frac{3\alpha_{1} A_{1} + \mu}{\alpha_{1}A_{1}}, \quad d =
\frac{(3\alpha_{1} - \rho) A_{1} + 2\mu}{\alpha_{1} A_{1}}\\[.2pc]
&A_{1} = \frac{1}{\rho -\alpha_{1}} \left[ \mu \left( 1 -
\frac{1}{K}\right) - H\right].
\end{align*}

In order to find $\eta (x)$ we make backward substitutions and finally
obtain the form
\begin{equation}
\eta (x) = \frac{A_{3}}{m_{1}(1 + m_{1})} {\rm e}^{(1 + m_{1}) x} +
\frac{A_{4}}{m_{2} (1 + m_{2})} {\rm e}^{(1 + m_{2}) x} + A_{5} {\rm
e}^{x} + A_{6}x,
\end{equation}
where $A_{i}\ (i = 5,6)$ are constants of integration and
\begin{equation*}
m_{1} = \frac{- c - \sqrt{c^{2} - 4d}}{2}, \quad m_{2} = \frac{- c +
\sqrt{c^{2} - 4d}}{2}.
\end{equation*}

From (3.2.4), (3.3.1) and (3.3.13) we get
\begin{align}
\psi (x, y) &= \frac{y}{\rho - \alpha_{1}}\left[ \mu \left( 1 -
\frac{1}{K} \right) - H\right] + A_{5}{\rm e}^{x} + A_{6}\nonumber\\[.2pc]
&\quad\ + \frac{A_{3}}{m_{1} (1 + m_{1})} {\rm e}^{(1 + m_{1}) x} +
\frac{A_{4}}{m_{2}(1 + m_{2})}{\rm e}^{(1 + m_{2}) x}.
\end{align}
The velocity components and pressure field are
\begin{align}
u &= \frac{1}{\rho - \alpha_{1}}\left[ \mu \left( 1 - \frac{1}{K}
\right) - H\right],\\[.2pc]
v &= -\left[ \frac{A_{3}}{m_{1}} {\rm e}^{(1 + m_{1}) x} +
\frac{A_{4}}{m_{2}} {\rm e}^{(1 + m_{2}) x} + A_{5}{\rm e}^{x}\right],\\[.2pc]
p &= p_{2} - \frac{1}{2}\rho [A_{1}^{2}]\nonumber\\[.2pc]
&\quad\  + \alpha_{1} \left[ \begin{array}{@{}c@{}}
\frac{A_{3}^{2}}{m_{1}^{2}} {\rm e}^{2(1 + m_{1}) x} +
\frac{2A_{3} A_{4}}{m_{1}m_{2}}{\rm e}^{(2 + m_{1} + m_{2}) x} +
\frac{A_{4}^{2}}{m_{2}^{2}} {\rm e}^{2 (1 + m_{2}) x} + A_{5}^{2} {\rm
e}^{2x}\\[.7pc]
+ \frac{2A_{3} A_{5} (2 + 3m_{1} + m_{1}^{2})}{(2 +
m_{1}) m_{1}} {\rm e}^{(2 + m_{1}) x} + \frac{2A_{4} A_{5}(2 +
3m_{2} + m_{2}^{2})}{(2 + m_{2}) m_{2}} {\rm e}^{(2 + m_{2}) x}
\end{array} \right],
\end{align}
where $p_{2}$ is the reference pressure.

\pagebreak

The streamline for $\psi = \Omega_{5}$ is given by the functional form
\begin{equation}
y = -\frac{1}{\varepsilon_{1}} \left[ \begin{array}{@{}c@{}}
- \Omega_{5} + \frac{A_{3}}{m_{1} (1 + m_{1})^{2}} {\rm e}^{(1 + m_{1})
x}\\[.7pc]
+ \frac{A_{4}}{m_{2} (1 + m_{2})^{2}} {\rm e}^{(1 + m_{2}) x} + A_{5}
{\rm e}^{x} + A_{6}x
\end{array} \right],
\end{equation}
where
\begin{equation*}
\varepsilon_{1} = \frac{1}{1 - \Lambda}\left[ \nu \left( 1 -
\frac{1}{K}\right) - \chi \right].
\end{equation*}
Streamline pattern is plotted in figure~7 for $\sigma =\lambda =1$, $\mu
/\rho = 0.5$, $\alpha_{1}/\rho = 0.1$, $K = 0.5$, $N = 0$, $\varphi =
0.05$, $A3 = A4 = A5 = A6 = 1,\ \psi = 15, 20, 25, 30, 40$.}
\end{case}

\begin{fig}
\hskip 4pc{\epsfxsize=8.7cm\epsfbox{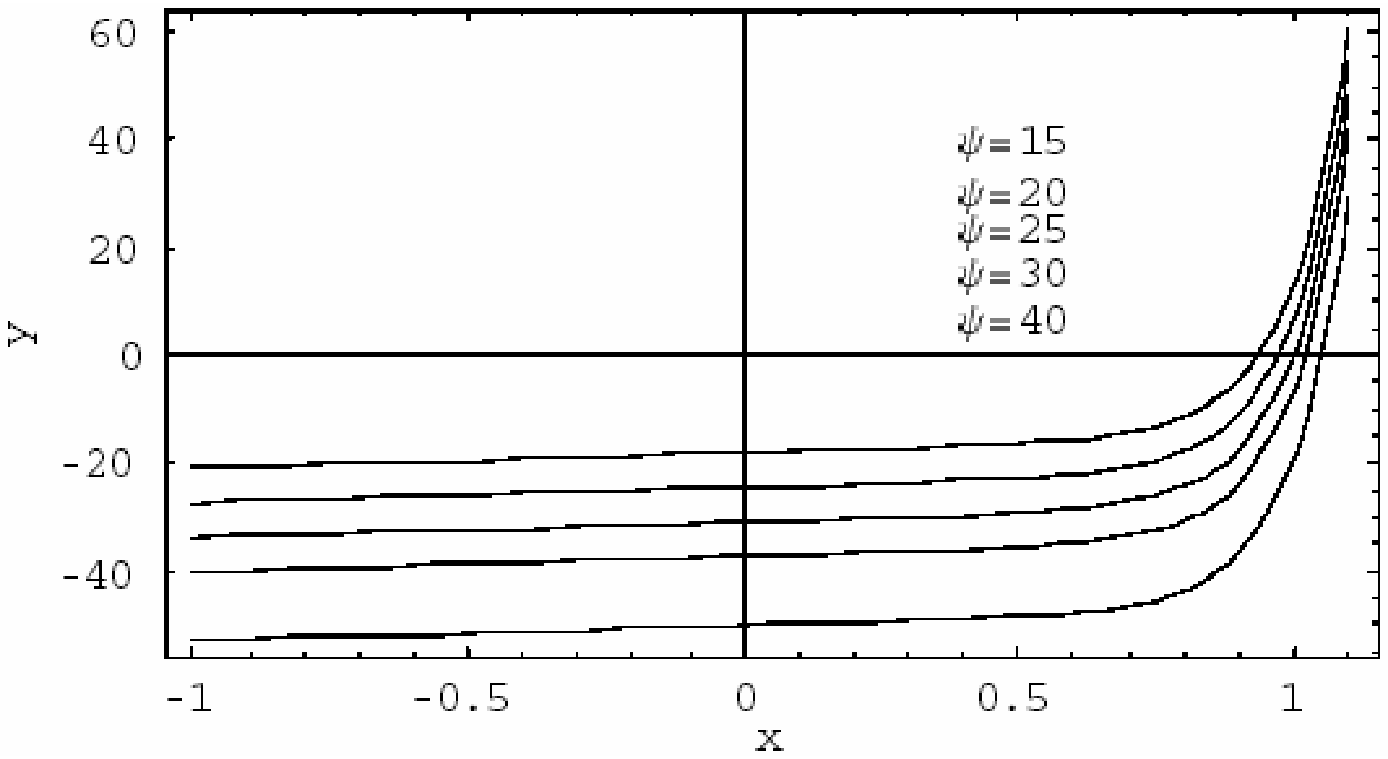}}\vspace{-.7pc}
\caption{Streamline flow pattern for $\protect\psi (x, y) =
\frac{y}{\protect\rho -\protect \alpha_{1}} \left[ \protect\mu \left( 1
- \frac{1}{K}\right) - H\right] + A_{5} {\rm e}^{x}$}
\hskip 4pc {\xxxx $  + A_{6} +
\frac{A_{3}}{m_{1}(1 + m_{1})} {\rm e}^{(1 + m_{1}) x} +
\frac{A_{4}}{m_{2}(1 + m_{2})} {\rm e}^{(1 + m_{2}) x}$.}\vspace{.9pc}
\end{fig}

\begin{case}{\rm
$\alpha_{1}\neq 0, \textbf{\ } \sigma =1, \textbf{\ } \lambda \neq 0$

Equation (3.3.7) reduces to
\begin{align}
&\alpha_{1} A_{1} (1 + \lambda {\rm e}^{x}) \eta^{\rm IV} + (\mu -
\alpha_{1} A_{1}\lambda {\rm e}^{x}) \eta^{\prime \prime \prime} + A_{1}
[(\alpha_{1} - \rho) \lambda {\rm e}^{x} - \rho] \eta^{\prime
\prime}\nonumber\\[.2pc]
&\quad\ = -\left[ (\rho -\alpha_{1}) A_{1}\lambda {\rm e}^{x} - \left(
H + \frac{\mu}{K}\right) \right] \eta^{\prime}.
\end{align}

To obtain the solution of (3.3.19) we try to reduce its order. For this
purpose we put $\eta^{\prime} = \widehat{A} (x)$ which leaves (3.3.19)
into a form which is one order less, that is
\begin{align}
&\alpha_{1} (1 + \lambda {\rm e}^{x}) \widehat{A}^{\prime \prime \prime}
+ (\mu/A_{1} - \alpha_{1}\lambda {\rm e}^{x}) \widehat{A}^{\prime
\prime} + [(\alpha_{1} - \rho) \lambda {\rm e}^{x} - \rho]
\widehat{A}^{\,\prime}\nonumber\\[.2pc]
&\quad\ = -\left[ (\rho -\alpha_{1}) \lambda {\rm e}^{x} - \frac{1}{A_{1}}
\left( H + \frac{\mu}{K}\right) \right] \widehat{A} = 0.
\end{align}
Now substituting $\widehat{A} (x) = \overline{P} (x) {\rm e}^{x}$ in
(3.3.20) and then $\overline{P}^{\prime} (x) = R(x)$ into the resulting
expression, we get
\begin{align}
&\alpha_{1} (1 + \lambda {\rm e}^{x}) R^{\prime \prime} + \left[ \frac{
K\mu (\rho -\alpha_{1})}{K(\mu - H) -\mu} + \alpha_{1} (3 + 2\lambda
{\rm e}^{x}) \right] R^{\prime}\nonumber\\[.2pc]
&\quad\ = -\left[ \frac{2K\mu (\rho -\alpha_{1})}{K(\mu - H) -\mu} +
(2\alpha_{1} - \rho) (1 + \lambda {\rm e}^{x}) \right] R,
\end{align}
where we have taken the constant of integration equal to zero.}
\end{case}

\begin{subcase}{\rm
The solution of (3.3.21) for $\lambda = 0$\ is given by
\begin{align}
R (x) &= C_{1}\exp \left( \frac{-X_{1} - \sqrt{X_{1}^{2} - 4X_{2}}}{2}
\right) x\nonumber\\[.2pc]
&\quad\ + C_{2}\exp \left( \frac{-X_{1} + \sqrt{X_{1}^{2} -
4X_{2}}}{2}\right) x,
\end{align}
where $C_{1}$ and $C_{2}$ are arbitrary constants and
\begin{equation}
X_{1} = \frac{K\mu (\rho -\alpha_{1})}{\alpha_{1}[K (\mu - H) - \mu]} +
3, \text{ } X_{2} = \frac{2K\mu (\rho - \alpha_{1})}{\alpha_{1} [K(\mu -
H) -\mu]} + \frac{3\alpha_{1} - \rho}{\alpha_{1}}.
\end{equation}
The backward substitution gives the value of $\eta (x)$ as
\begin{equation}
\eta (x) = \frac{C_{1}}{\overline{m}_{1} (1 + \overline{m}_{1})} {\rm
e}^{(1 + \overline{m}_{1}) x} + \frac{C_{2}}{\overline{m}_{2} (1 +
\overline{m}_{2})} {\rm e}^{(1 + \overline{m}_{2}) x} + C_{3} {\rm
e}^{x} + C_{4},
\end{equation}
where $C_{i}\ (i = 1, 2, 3, 4)$ are constants of integration and
\begin{equation*}
\overline{m}_{1} = \frac{-X_{1} - \sqrt{X_{1}^{2} - 4X_{2}}}{2},\quad
\overline{m}_{2} = \frac{-X_{1} + \sqrt{X_{1}^{2} - 4X_{2}}}{2}.
\end{equation*}

The stream function, the velocity components and the pressure field in
this case are respectively given as
\begin{align}
&\psi (x, y) = \frac{y}{\rho -\alpha_{1}}\left[ \mu \left( 1 -
\frac{1}{K} \right) - H\right] + C_{3} {\rm e}^{x} + C_{4}\nonumber\\[.2pc]
&\hskip 3.6pc + \frac{C_{1}}{\overline{m}_{1} (1 + \overline{m}_{1})} {\rm
e}^{(1 + \overline{m}_{1}) x} + \frac{C_{2}}{\overline{m}_{2} (1 +
\overline{m}_{2})} {\rm e}^{(1 + \overline{m}_{2}) x},\\[.2pc]
&u = \frac{1}{\rho -\alpha_{1}}\left[ \mu \left( 1 - \frac{1}{K}\right)
- H\right],\\[.2pc]
&v = -\left[ \frac{C_{1}}{\overline{m}_{1}} {\rm e}^{(1 +
\overline{m}_{1}) x} + \frac{C_{2}}{\overline{m}_{2}} {\rm e}^{(1 +
\overline{m}_{2}) x} + C_{3} {\rm e}^{x}\right],\\[.2pc]
&p = p_{3} - \frac{1}{2}\rho [A_{1}^{2}]\nonumber\\[.2pc]
&\hskip 1.4pc +\!\alpha_{1}\! \left[ \!\begin{array}{@{}c@{}}
\frac{C_{1}^{2}}{\overline{m}_{1}^{2}} {\rm e}^{2(1 +
\overline{m}_{1})x} + \frac{2C_{1}C_{2}}{\overline{m}_{1} \overline{m}_{2}}
{\rm e}^{(2 + \overline{m}_{1} + \overline{m}_{2}) x} +
\frac{C_{2}^{2}}{\overline{m}_{2}^{2}} {\rm e}^{2(1 + \overline{m}_{2})
x} + C_{3}^{2}{\rm e}^{2x}\\[.7pc]
+ \frac{2C_{1}C_{3}(2 + 3\overline{m}_{1} +
\overline{m}_{1}^{2})}{(2 + \overline{m}_{1}) \overline{m}_{1}}{\rm
e}^{(2 + \overline{m}_{1}) x} + \frac{2C_{2}C_{3} (2 +
3\overline{m}_{2} + \overline{m}_{2}^{2})}{(2 + \overline{m}_{2})
\overline{m}_{2}} {\rm e}^{(2 + \overline{m}_{2}) x}
\end{array}\! \right],
\end{align}
where $p_{3}$ is the reference pressure and the streamline for $\psi =
\Omega_{6}$ is given by the functional form
\begin{equation}
y = -\frac{1}{\varepsilon_{1}} \left[ \begin{array}{@{}c@{}}
- \Omega_{6} + \frac{C_{1}}{\overline{m}_{1} (1 + \overline{m}_{1})^{2}}
{\rm e}^{(1 + \overline{m}_{1}) x}\\[.7pc]
+ \frac{C_{2}}{\overline{m}_{2} (1 + \overline{m}_{2})^{2}}
{\rm e}^{(1 + \overline{m}_{2}) x} + C_{3} {\rm e}^{x} + C_{4}
\end{array} \right].
\end{equation}
Streamline pattern is plotted in figure~8 for $\sigma = \lambda = 1$,
$\mu /\rho = 0.5$, $\alpha_{1}/\rho = 0.1$, $K = 1.8$, $N = 1$, $\varphi
= 1$, $C_{1} = C_{2} = C_{3} = C_{4} = 1,\ \psi = 15, 20, 25, 30, 40$.}
\end{subcase}

\begin{fig}
\hskip 4pc{\epsfxsize=8.7cm\epsfbox{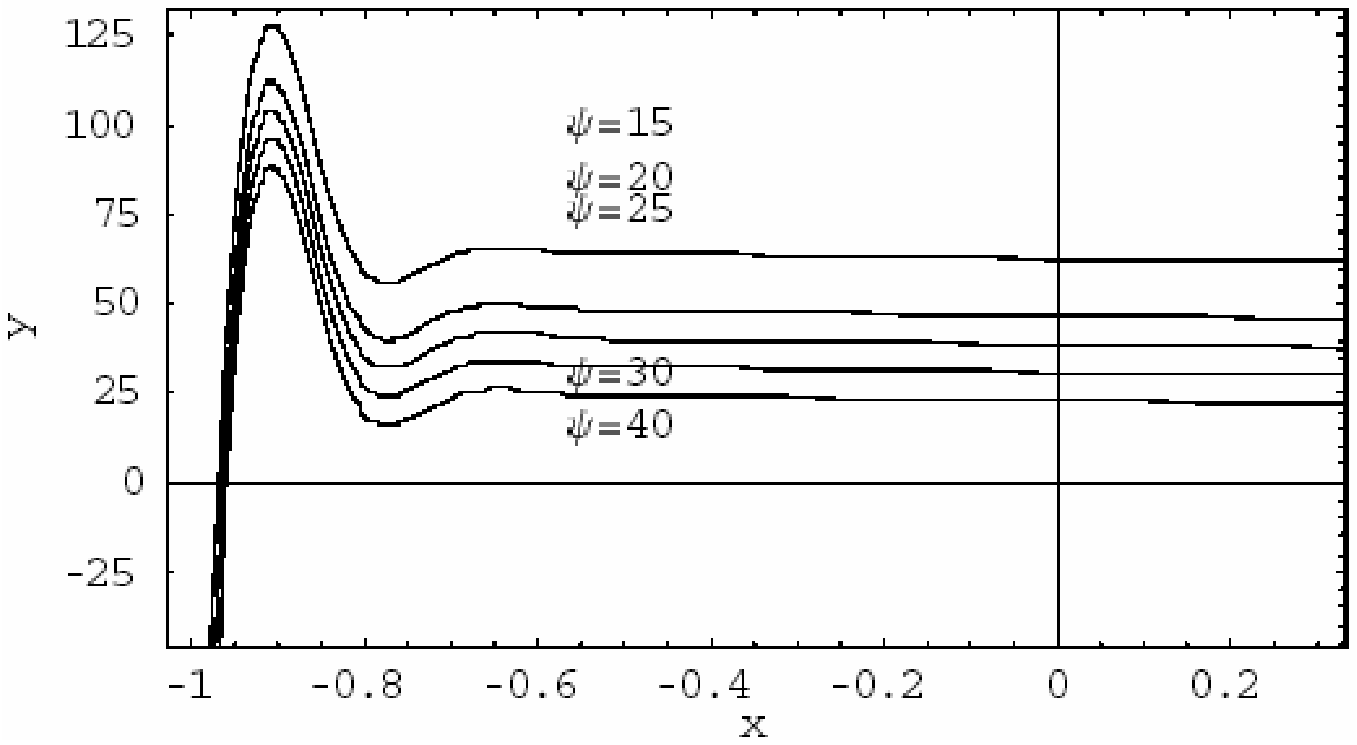}}\vspace{-.7pc}
\caption{Streamline flow pattern for $\protect\psi (x, y) =
\frac{y}{\protect\rho -\protect \alpha_{1}}[\protect\mu (1 -
\frac{1}{K}) - H] + C_{3} {\rm e}^{x} + $}
\hskip 4pc {\xxxx $C_{4} +
\frac{C_{1}}{\overline{m}_{1} (1 + \overline{m}_{1})} {\rm e}^{(1 +
\overline{m}_{1}) x} + \frac{C_{2}}{\overline{m}_{2}(1 +
\overline{m}_{2})} {\rm e}^{(1 + \overline{m}_{2}) x}$.}\vspace{1pc}
\end{fig}

\begin{subcase}\hskip -.4pc {\rm (General case).\ \ We now try to obtain
the solution of (3.3.21) for $\lambda \neq 0$. Equation (3.3.21) may
be written as
\begin{equation}
(1 + \lambda {\rm e}^{x}) R^{\prime \prime} + [ X_{1} + 2\lambda {\rm
e}^{x}] R^{\prime} + [X_{2} + X_{3}\lambda {\rm e}^{x}] R = 0.
\end{equation}
In order to solve (3.3.30) we put $\theta = {\rm e}^{x}$ to get the
following equation
\begin{equation}
(1 + \lambda \theta) \theta^{2}R^{\prime \prime} + (1 + X_{1} + 2\lambda
\theta) \theta R^{\prime} + (X_{2} + X_{3}\lambda \theta) R = 0,
\end{equation}
where differentiation is with respect to $\theta$. The solution for
(3.3.31) is obtained through {\it Mathematica} and is given by
\begin{align*}
R (\theta) &= \theta^{1/2} \left(-X_{1} - \sqrt{X_{1}^{2} - 4X_{2}}\right)\\[.3pc]
&\quad\ \times \left[ \begin{array}{@{}c@{}}
C_{5} \theta^{\sqrt{X_{1}^{2} - 4X_{2}}} 2F1 \left\{ \frac{1}{2}
\Phi_{1}, \frac{1}{2} \Phi_{2}, 1 + \sqrt{X_{1}^{2} - 4X_{2}}, - \theta
\lambda \right\}\\[1pc]
+ C_{1} 2F1 \left\{ \frac{1}{2} \Phi_{3}, \frac{1}{2} \Phi_{4}, 1 -
\sqrt{X_{1}^{2} - 4X_{2}}, - \theta \lambda \right\}
\end{array} \right],
\end{align*}
where
\begin{align*}
\Phi_{1} &= 2-X_{1} + \sqrt{X_{1}^{2} - 4X_{2}} - 2\sqrt{1-X_{3}},\\[.2pc]
\Phi_{2} &= 2-X_{1} + \sqrt{X_{1}^{2} - 4X_{2}} + 2\sqrt{1-X_{3}},\\[.2pc]
\Phi_{3} &= 2-X_{1} - \sqrt{X_{1}^{2} - 4X_{2}} - 2\sqrt{1-X_{3}},\\[.2pc]
\Phi_{4} &= 2-X_{1} - \sqrt{X_{1}^{2} - 4X_{2}} + 2\sqrt{1-X_{3}},
\end{align*}
and $_{2}F_{1}$ is the hypergeometric function defined in Appendix~A. The
stream function, velocity components and the pressure field can be obtained
through the definitions of $_{2}F_{1} [a, b, c, z]$.}
\end{subcase}

\section{Concluding remarks}

In this paper, the exact solutions of non-linear equations governing the
flow for a second-grade fluid in a porous medium are obtained by
assuming different forms of the stream function (already used by various
authors in different situations), in presence of a strong magnetic
field. The expressions for velocity profile, streamline and pressure
distribution are constructed in each case. Our result indicates that
velocity, stream function and pressure are strongly dependent upon the
material parameter $\alpha_{1}$ of the second-grade fluid. It is shown
through graphs that increase in $\alpha_{1}$ leads to decrease in
velocity and decrease in $\alpha_{1}$ leads to increase in velocity (see
figures~4 and 5). Also, the present analysis is more general and several
results of various authors (as already mentioned in the text) can be
recovered in the limiting cases.

\section*{Appendix A}

Hypergeometric $_{2}F_{1}[a, b, c, z]$ is the hypergeometric function
$_{2}F_{1} [a, b; c; z]$ and is the special case of the generalized
hypergeometric function $_{p}F_{q} [\mathbf{a}; \mathbf{b}; z]$ for
$p=1$ and $q=1$.

Hypergeometric function has the following properties:
\begin{enumerate}
\renewcommand\labelenumi{\arabic{enumi}.}
\item The $_{2}F_{1}$ function has the series expansion $_{2}F_{1} [a,
b; c; z] = \sum_{k=0}^{\infty}
\frac{(a)_{k}(b)_{k}}{(c)_{k}}\frac{z^{k}}{k!}$.

\item Hypergeometric $_{2}F_{1} [a, b, c, z]$ has a branch cut
discontinuity in the complex $z$-plane running from $1$ to $\infty$.
\end{enumerate}

\end{document}